\documentclass[11pt,a4]{article}

\usepackage{amsmath,amssymb,amsthm}
\usepackage{color}
\usepackage{graphicx}

\newtheorem{definition}{Definition}

\newtheorem{proposition}{Proposition}
\newtheorem{remark}{Remark}

\def \beq{ \begin{equation} }
\def \eeq{\end{equation}}

\newcommand{\e}{\epsilon}
\title{Three-body relative equilibria on $S^2$\\
II: Extended Lagrangian configurations}
\date{}  
\begin{document}
	\maketitle
	\author{\begin{center}
	{ Toshiaki~Fujiwara$^1$, Ernesto P\'{e}rez-Chavela$^2$}\\	
		\bigskip
	   $^1$College of Liberal Arts and Sciences, Kitasato University,       Japan. fujiwara@kitasato-u.ac.jp\\
	    $^2$Department of Mathematics, ITAM, M\'exico. ernesto.perez@itam.mx\\
	\end{center}
	


\begin{flushright}
{\it To the memory of our friend Florin Diacu,\\
who was the inspiration of this work}
\end{flushright}

\begin{abstract}
This is a natural continuation of our first paper \cite{pre}, where we develop a new geometrical technique which allow us to study relative equilibria on the two sphere. We consider a system of three positive masses on $\mathbb{S}^2$  moving 
under the influence of an generic attractive potential which only depends on the mutual distances among the masses. 
We reduce the problem of finding extended Lagrangian relative equilibria to the analysis of the inertia tensor, then we obtain a more manageable equivalent inertia tensor which allow us to find new families of Lagrangian configurations.
\end{abstract}


\section{Introduction}
The simplest motions in dynamics of particles are those where the mutual distances among the particles remain constant for all time, like if they 
constituted a rigid body, these motions are called relative equilibria, that for short we denote simply as $RE$. 
In a previous work \cite{pre}, we have studied relative equilibria, for a system of three positive masses moving on the two dimensional sphere
$\mathbb{S}^2$, under the influence of a general attractive potential which only depends on the mutual distances among the masses, in that paper we analyze the case where the masses are on the same geodesic (a great circle of   $\mathbb{S}^2$), we call them Euler relative equilibria.

In this paper, we tackled the relative equilibria where the masses are not on the same geodesic, we call them {\it extended Lagrangian relative equilibria} or {Lagrangian configurations} by short, ahead in this paper we explain this name. The main achievement of this article is had introduced a new geometric technique to study them. 
In this way it is important to distinguish 
two concepts ``shape'' and `configuration''
to investigate relative equilibria on $\mathbb{S}^2$.
The shape stands for the set of
mutual angles $\sigma_{ij}$ 
between the bodies on $\mathbb{S}^2$.
The configuration stands for the set of coordinates
of three bodies $\theta_k$ and $\phi_k$ 
(in a spherical coordinate system). 
The map from the configuration to the shape is trivial.
However, the inverse, namely, 
the map from the shape to the configuration is not trivial.
The new geometrical technique consists in giving a 
map from the shape to the configuration,
what we call the ``translation formula''.
Utilising this formula,
we first derive some conditions  
for a shape to satisfy the equations of motion.
If the shape satisfies the conditions,
it is called a ``rigid rotator''.
Then, utilising the translation formula again,
we get the configuration for the rigid rotator,
namely the relative equilibrium configuration.

We can apply this technique to general attractive potentials which only depend on the mutual distances among the masses. The cotangent potential, that is generally used to study relative equilibria on spaces of constant positive curvature is included in this family of potentials, we use it to exemplify some concrete applications of our technique. In this way we can compare some similar results obtained by other authors using different approaches.
Until we know, just for the cotangent potential, the first example of this kind of $RE$ appears in \cite{Diacu-EPC1}, where the authors proved the existence of $RE$ for three equal masses located at the vertices of an equilateral triangle moving on a plane parallel to the equatorial plane. After that, some other authors have also studied Lagrangian configurations, but just for equilateral triangle shapes moving on a plane parallel to the equatorial plane in the context of the curved positive problem (that is using the cotangent potential), see for instance \cite{Diacu1, Diacu4, M-S, EPC1, tibboel}. In a recent work, the authors were able to obtain a continuation of any relative equilibria of the Newtonian $n$--body problem to spaces of constant curvature, when the value of the curvature is small \cite{Bengochea}. In particular they could extend Lagrangian $RE$ with arbitrary masses to the sphere.
Other authors, by using a totally different techniques have studied $RE$ on $\mathbb{S}^2$ for the three vortices problem, see for instance \cite{Borisov1} and the references therein.

By using our technique we extend the result in \cite{Diacu-EPC1} to any attractive potential depending only on the mutual distances among the masses. But this is just an easy application that shows how powerful our technique is. In general it will allow us to find many new families of Lagrange configurations,
for the case of three positive masses moving on $\mathbb{S}^2$,
results that are totally original. We will describe this technique in the following, when we break down the content of this article.

The paper is organized as follows: in Section \ref{notations}, we introduce the equations of motion for the three body problem on $\mathbb{S}^2$ in spherical coordinates, compute the angular momentum ${\bf c}=(c_x,c_y,c_z)$, where $c_x=c_y=0$, and define the Eulerian and the extended Lagrangian relative equilibria. We have used the same notations and equation of motion of our previous paper \cite{pre}, here we just give a summary of that material. Section \ref{eqLRE} is dedicated to compute the equations of motion for the extended Lagrangian $RE$, and to prove that, in order to have these motions, the three masses must be on the same hemisphere. Section \ref{I-tensor} contains the main results of this paper. We start defining the inertia tensor $I$, which is a symmetric real matrix, then the eigenvalue problem 
$$I\psi_{\alpha} = \lambda_{\alpha}\psi_{\alpha}$$
has three eigenvalues and three mutually orthogonal eigenvectors called principal axis. We prove that the $z$--axis is one of the principal axes, and from here we reduce the problem to find a rotation axes to calculate the inertia tensor $I$ for a given shape. This can be a really tedious problem, however we find an alternative inertia tensor $J$ which is similar to $I$. This is the new geometrical method that we have developed to study Lagrangian rigid rotators. In 
Section \ref{LRE}, we find some extended Lagrangian $RE$ for the case of three equal masses for generic potentials, then in order to find other concrete families of Lagrangian $RE$ we restrict our analysis of the cotangent potential.

\section{Notations and equations of motion}\label{notations}
We are using the same notations, definitions
and equations of motion that in our previous paper \cite{pre}. 
 In order to have a self contained paper we summarize in this section that material.

\subsection{Notations for the spherical coordinates}
To study the dynamics on $\mathbb{S}^2$ we use 
spherical coordinate $(\theta,\phi)$.
The relation  between the Cartesian coordinates $(X,Y,Z)$
and the spherical coordinates is
\begin{equation}
(X,Y,Z)=R(\sin\theta\cos\phi,\sin\theta\sin\phi,\cos\theta).
\end{equation}
The chord length $D_{ij}$ and 
the arc angle $\sigma_{ij}$ between the  points
$(X_i,Y_i,Z_i)$ and $(X_j,Y_j,Z_j)$ is given by
\begin{equation}
\label{defD}
D_{ij}^2
=2R^2\Big(1-\cos\theta_i\cos\theta_j-\sin\theta_i\sin\theta_j\cos(\phi_i-\phi_j)
		\Big),
\end{equation}
\begin{equation}
\label{relationDandSigma}
\e D_{ij}=\sin(\sigma_{ij}/2) \quad
\mbox{ with } \quad
\e =1/(2R).
\end{equation}
The arc angle is equal to the angle between the two points 
as seen from the origin of $\mathbb{S}^2$, therefore we have
\begin{equation}
0\le \sigma_{ij}\le \pi.
\end{equation}

The equations (\ref{defD}) and (\ref{relationDandSigma}) yield
the fundamental relation
for the arc angle and the spherical coordinates.
\begin{equation}
\label{fundamentalrelation}
\cos\sigma_{ij}
=\cos\theta_i\cos\theta_j+\sin\theta_i\sin\theta_j\cos(\phi_i-\phi_j).
\end{equation}

\subsection{Equations of motion}
The Lagrangian for the three body problem on
$\mathbb{S}^2$
is given by
\begin{equation}
\label{theLagrangian}
\begin{split}
&L=K+V, \quad \text{where}\\
&K=R^2\sum_{k=1,2,3}\frac{m_k}{2}
		\left(\dot{\theta}_k^2+\sin^2(\theta)\dot{\phi}_k^2\right),\quad
V=\sum_{i<j}m_i m_j U(D_{ij}^2),\\
\end{split}
\end{equation}
here dot on symbols represents 
time derivative.
Along the paper we will use the  notation for the derivative
$U'$ which is defined by
\begin{equation}
U'(D^2)=\frac{dU(D^2)}{d(D^2)}<0.
\end{equation}
The last inequality is the assumption that we make
to ensure that the force between any two bodies is attractive.

When we need to specify a particular potential,
we use the cotangent potential and its derivative,
\begin{equation}
U(D_{ij}^2)
=\frac{1-2\e^2D_{ij}^2}{\sqrt{D_{ij}^2(1-\e^2D_{ij}^2)}}
=\frac{\cos\sigma_{ij}}{R\sqrt{1-\cos^2(\sigma_{ij})}},
\end{equation}
\begin{equation}\label{deriv-pot}
U'(D^2)=-\frac{1}{2R^3\sin^3(\sigma)}.
\end{equation}

Since the Lagrangian is  invariant under $SO(3)$ rotation
around the centre of $\mathbb{S}^2$,
the angular momentum 
\begin{equation}
{\bf c} = (c_x, c_y, c_z)=\sum_k m_k (X_k,Y_k,Z_k)\times (\dot{X}_k, \dot{Y}_k, \dot{Z}_k)
\end{equation}
is a first integral.
Each component is represented as
\begin{align}
c_x&=R^2 \sum_k m_k \left(-\sin(\phi_k) \dot\theta_k
	-\sin(\theta_k)\cos(\theta_k)\cos(\phi_k) \dot\phi_k\right),
	\label{defCx}\\
c_y&=R^2 \sum_k m_k \left(\cos(\phi_k) \dot\theta_k
	-\sin(\theta_k)\cos(\theta_k)\sin(\phi_k) \dot\phi_k\right),
	\label{defCy}\\
c_z&=R^2 \sum_k m_k \sin^2(\theta_k) \dot \phi_k.
	\label{defCz}
\end{align}

\subsection{Relative equilibria and the equations of motion for them}
\begin{definition}[Relative equilibrium]
An relative equilibrium on $\mathbb{S}^2$ is a solution 
of the equations of motion that satisfy
$\dot{\theta}_k=0$ and $\dot{\phi}_i-\dot{\phi}_j=0$
for all $k=1,2,3$ and all pair $(i,j)$
in a spherical coordinate system, we will call them as $RE$.
\end{definition}

By the above definition of $RE$, 
the angular momentum has the form,
\begin{equation}
{\bf c} = (0,0,c_z),
\end{equation}
\begin{equation}
c_z=R^2\omega\sum_k m_k \sin\theta_k\cos\theta_k,
\end{equation}
$c_x=c_y=0$ are expressed as
\begin{equation}
\label{CxCy}
\sum_k m_k\sin(\theta_k)\cos(\theta_k) e_k=0
\mbox{ if } \omega\ne 0,
\end{equation}
where $e_k$ is the two-dimensional unit vector
\begin{equation}
e_k=(\cos\phi_k,\sin\phi_k).
\end{equation}

When $\omega=0$, the corresponding relative equilibrium is a fixed point.

Taking the inner product between equation \eqref{CxCy} and $e_k$ we obtain
\begin{equation}
\label{cxcyInE3}
\begin{split}
m_i\sin(\theta_i)\cos(\theta_i)\cos(\phi_k-\phi_i)
+m_j\sin(\theta_j)\cos(\theta_j)\cos(\phi_k-\phi_j)\\
+m_k\sin(\theta_k)\cos(\theta_k)
=0,
\end{split}
\end{equation}
for $(i,j,k)=(1,2,3),(2,3,1)$, and $(3,1,2)$.
This cyclic notation will be used
when the three indexes $(i,j,k)$ appear in one equation.
Now, taking the outer product between the same equation and $e_k$ we get
\begin{equation}
\label{cxcyInE4}
m_i\sin(\theta_i)\cos(\theta_i)\sin(\phi_k-\phi_i)
+m_j\sin(\theta_j)\cos(\theta_j)\sin(\phi_k-\phi_j)
=0.
\end{equation}

The equations of motion for $\phi_k$ and $\theta_k$ 
for the relative equilibria are reduced to
\begin{equation}
\label{relationsFromPhiII}
\begin{split}
&m_1 m_2 U'(D_{12})\sin\theta_1 \sin\theta_2 \sin(\phi_1-\phi_2)\\
=&m_2 m_3 U'(D_{23})\sin\theta_2 \sin\theta_3 \sin(\phi_2-\phi_3)\\
=&m_3 m_1 U'(D_{31})\sin\theta_3 \sin\theta_1 \sin(\phi_3-\phi_1)
\end{split}
\end{equation}
and
\begin{equation}
\label{eqForThetaII}
\begin{split}
-&\omega^2 m_k \sin\theta_k\cos\theta_k\\
&=2\sum_{i\ne k}m_k m_iU'(D_{ki}^2)
	\Big(
		\sin\theta_k\cos\theta_i-\cos\theta_k\sin\theta_i\cos(\phi_k-\phi_i)
	\Big).
\end{split}
\end{equation}

The Eulerian and Extended Lagrangian relative equilibria
are defined as follows.
\begin{definition}[Eulerian and Extended Lagrangian relative equilibria]\label{defRE}
An ``Eulerian relative equilibrium'' is
a relative equilibrium in which the three bodies are on
the same geodesic (a great circle) of $\mathbb{S}^2$.
The other equilibria are called
``Extended Lagrangian relative equilibria'' or simply $LRE$ by short.
\end{definition}


\section{Equations of motion for the extended Lagrangian relative
equilibria}\label{eqLRE}

First we observe that from equation (\ref{relationsFromPhiII}), if
$\sin\theta_i\sin\theta_j\sin(\phi_i-\phi_j)=0$ for one pair $(i,j)$,
then all of them must be zero since
 $U'(D_{ij}^2)\ne 0$ for all $i,j$. 
Then, in this case the three bodies must be on a same meridian. To verify this statement we observe that if no bodies are on the poles, then $\sin \theta_k \ne 0$ for all $k$, and the result follows easily. 
If just one body is at one of the poles, say $\sin \theta_3=0$ and $\sin \theta_1 \ne 0$, $\sin \theta_2 \ne 0$, then $\sin (\phi_1-\phi_2) = 0$, from here 
$\phi_1 - \phi_2=0$ or $\pi$.
and therefore the three bodies are on a same meridian since the poles belong to any meridian.
If at least two bodies are on the poles, the result follows immediately. 

Therefore, by Definition \ref{defRE},
the extended Lagrangian relative equilibria must satisfy:
\begin{equation}
\label{conditionForLagrangian1}
\sin\theta_k\ne 0
\mbox{ and }
\sin(\phi_i-\phi_j)\ne 0
\mbox{ for all }k
\mbox{ and all pair }i,j.
\end{equation}
Then,
by the definition of spherical coordinates $\sin\theta_k>0$,
and by the property of the potential $U'(D_{ij}^2)<0$,
all $\sin(\phi_i-\phi_j)$ must have the same sign.
Without loss of generality, we assume
\begin{equation}
\sin(\phi_i-\phi_j)>0 \quad
\textcolor{red} \quad 
\mbox{ for } \quad (i,j)=(1,2), (2,3), (3,1).
\end{equation}

 In the following calculations, we will use the equation
 (\ref{cxcyInE4}) with $(i,j,k)=(3,1,2)$, that is
 \begin{equation}
 \label{cxcyInE5}
 m_3\sin\theta_3\cos\theta_3\sin(\phi_2-\phi_3)
 =m_1\sin\theta_1\cos\theta_1\sin(\phi_1-\phi_2).
 \end{equation}
Multiplying by $\cos(\theta_3)$ both sides of
the first and the second expressions in equation (\ref{relationsFromPhiII}),
we obtain
\begin{equation}
\begin{split}
&m_1m_2U'(D_{12}^2)\sin\theta_1\sin\theta_2\sin(\phi_1-\phi_2)\cos\theta_3\\
&=m_2U'(D_{23}^2)\sin\theta_2
(m_3\sin\theta_3\cos\theta_3\sin(\phi_2-\phi_3))\\
&=m_2U'(D_{23}^2)\sin\theta_2
	(m_1\sin\theta_1\cos\theta_1\sin(\phi_1-\phi_2)).
\end{split}
\end{equation}
To get the last line, we have used (\ref{cxcyInE5}).
Dividing the first and the last line by
$m_1m_2\sin\theta_1\sin\theta_2\sin(\phi_1-\phi_2)\ne 0$,
we obtain $U'(D_{12}^2)\cos\theta_3=U'(D_{23}^2)\cos\theta_1$.
Similarly, we  obtain,
\begin{equation}
\label{conditionForLagrangian2}
U'(D_{12}^2)\cos\theta_3
=U'(D_{23}^2)\cos\theta_1
=U'(D_{31}^2)\cos\theta_2.
\end{equation}

Again, if one $\cos\theta_k=0$, then all $\cos\theta_k$ must be zero since $U'(D_{ij}^2)\ne 0$.
This configuration is a collinear configuration on the Equator.
Therefore, for the extended Lagrangian equilibria
\begin{equation}
\label{conditionForLagrangian3}
\cos\theta_k\ne 0.
\end{equation}

Thus the conditions for the extended Lagrangian
relative equilibria are (\ref{conditionForLagrangian1}),
(\ref{conditionForLagrangian2}), and
(\ref{conditionForLagrangian3}).

Before proceeding further, we explain why the solution in 
(\ref{conditionForLagrangian2})
deserves the name ``Lagrangian''.
Consider the Euclidean limit near the North Pole,
$R\to \infty$ with $r_k=R\theta_k$ finite.
Then $\cos\theta_k=1+O(R^{-2})$.
Therefore for this limit,
$U'(D_{12}^2) \simeq U'(D_{23}^2) \simeq U'(D_{31}^2)$,
namely $\sigma_{12}\simeq \sigma_{23}\simeq \sigma_{31}$,
that is the equilateral triangle neglecting $O(R^{-2})$ terms.
So, it deserves the name ``Lagrangian''.

Therefore, the solutions of (\ref{conditionForLagrangian2})
contain nearly equilateral triangle for sufficiently small $\theta_k$.
However, the solutions of (\ref{conditionForLagrangian2})
can contain much more non-trivial shapes.

Now, we proceed further from (\ref{conditionForLagrangian2}).
Since, $U'(D_{ij}^2)<0$ for any $D_{ij}$,
the three quantities $\cos\theta_k$ must have the same sign, then we obtain the following result.
\begin{proposition}
\label{cosMustHaveTheSameSign}
The three bodies for any extended Lagrangian relative equilibrium
are on the same hemisphere.
\end{proposition}

Without loss of generality, we can assume
\begin{equation}
\cos\theta_k>0.
\end{equation}

Since $\cos\theta_k\ne 0$, dividing equation (\ref{conditionForLagrangian2})
by $\cos\theta_1\cos\theta_2\cos\theta_3$ we observe that the expression 
$U'(D_{ij}^2)/(\cos\theta_i\cos\theta_j)$ has a common value, that we denote as  $- \gamma$. That is
\begin{equation}
\label{conditionForLagrangian4}
U'(D_{ij}^2)=- \gamma\,\,\cos\theta_i\cos\theta_j \quad \text{for all} \quad (i,j).
\end{equation}

Putting the expressions 
(\ref{conditionForLagrangian4})
into the equation for $\theta_1$ in (\ref{eqForThetaII}) and using  (\ref{cxcyInE3}) we obtain
\begin{equation}
\label{equationForTheta1DetermineGamma}
\begin{split}
-\omega^2 m_1 &\sin\theta_1\cos\theta_1\\
=&-2\gamma\,m_1\sin\theta_1\cos\theta_1
	\Big(m_2\cos^2(\theta_2)+m_3\cos^2(\theta_3)\Big)\\
	&+2\gamma\, m_1 \cos^2(\theta_1)
		\Big(
			m_2\sin\theta_2\cos\theta_2\cos(\phi_1-\phi_2)
			+m_3\sin\theta_3\cos\theta_3\cos(\phi_3-\phi_1)
		\Big)\\
=&-2\gamma \,m_1\sin\theta_1\cos\theta_1
	\Big(m_2\cos^2(\theta_2)+m_3\cos^2(\theta_3)\Big)\\
	&+2\gamma \, m_1 \cos^2(\theta_1)
		\Big(
			-m_1\sin\theta_1\cos\theta_1
		\Big)\\
=&-2\gamma\,m_1\sin\theta_1\cos\theta_1
	\sum_k m_k \cos^2(\theta_k).
\end{split}
\end{equation}

Therefore, the value $\gamma$ is determined to be
\begin{equation}
\label{defOfC}
\gamma =\frac{\omega^2}{2\sum_k m_k \cos^2(\theta_k)}
>0.
\end{equation}

Similarly, we can easily check that the equations of motion
for $\theta_2$ and $\theta_3$ are also satisfied by the same $\gamma$.
Thus, all equations of motion are satisfied.

So far,
it is shown that
if the shape defined by $\sigma_{ij}$ satisfies
the equation (\ref{conditionForLagrangian4}) with
(\ref{defOfC}), it satisfies the equations of motion.
Namely, the shape is a rigid rotator.

The remaining problem is 
what is the correspondence
between $\sigma_{ij}$ and $\theta_k$.
In other words, how we can find the $z$--axis
for given shape. In the following section, by using equation (\ref{CxCy})
we will explain how do it. 

Before closing this section,
we would like to add a remark for a repulsive force.

\begin{remark}
We observe that the above discussion up to the equation (\ref{equationForTheta1DetermineGamma})
remains completely unchanged 
even for the repulsive force $U'(D^2)>0$,
because we have actually used the fact
$U'(D^2)$ has definite sign.
Then the equation of motion
demands that $\gamma$ must be positive
as shown in (\ref{defOfC}).
This means that $U'(D^2)$ must be negative
by (\ref{conditionForLagrangian4}).
Namely,
any repulsive force $U'(D^2)>0$
has no Lagrangian $RE$ on $\mathbb{S}^2$. 
\end{remark}

\section{The inertia tensor}\label{I-tensor}
The inertia tensor has an important role in the study of the rigid body problem as you can see in \cite{LandauLifshitz, Goldstein, Hestenes}.
Since for $RE$ the masses behave as a rigid body, the inertia tensor is also very important in the analysis of the RE on $\mathbb{S}^2$. In this section we first give the formal definition of this concept;  then we explain how we determine the axis of rotation in a 
$RE$. We finish this section getting a classification of the $RE$, we also show the correspondence between the shape variables $\sigma_{ij}$
and the angles $\theta_k$ for the extended Lagrangian relative equilibria.

\subsection{Definition of the inertia tensor}
The inertia tensor is defined by
\begin{equation}
I=\left(\begin{array}{ccc}
I_{xx} & I_{xy} & I_{xz} \\
I_{yx} & I_{yy} & I_{yz} \\
I_{zx} & I_{zy} & I_{zz}
\end{array}\right).
\end{equation}
Where 
\begin{equation}
\label{inertiaTensorDiagonal}
\begin{split}
I_{xx}&=\sum_{k} m_k (Y_k^2+Z_k^2)
	=\sum_k m_k \Big(\cos^2(\theta_k)+\sin^2(\theta_k)\sin^2(\phi_k)\Big),\\
I_{yy}&=\sum_{k} m_k (Z_k^2+X_k^2)
	=\sum_k m_k\Big(\cos^2(\theta_k)+\sin^2(\theta_k)\cos^2(\phi_k)\Big),\\
I_{zz}&=\sum_{k} m_k (X_k^2+Y_k^2)
	=\sum_k m_k \sin^2(\theta_k),\\
\end{split}
\end{equation}
and
\begin{equation}
\label{inertiaTensorOffDiagonal}
\begin{split}
I_{xy}=I_{yx}&= -\sum_{k} m_k X_k Y_k
	=-\sum m_k \sin^2(\theta_k)\sin\phi_k\cos\phi_k,\\
I_{xz}=I_{zx}&= -\sum_{k} m_k X_k Z_k
	=-\sum m_k \sin\theta_k\cos\theta_k\cos\phi_k,\\
I_{yz}=I_{zy}&= -\sum_{k} m_k Y_k Z_k
	=-\sum m_k \sin\theta_k\cos\theta_k\sin\phi_k.\\
\end{split}
\end{equation}

Since the inertia tensor is a symmetric real matrix,
the eigenvalue problem
\begin{equation}\label{eigenvalue-problem}
I\psi_\alpha=\lambda_\alpha\psi_\alpha
\end{equation}
has three eigenvalues 
and three mutually orthogonal eigenvectors.
These eigenvectors are called principal axes.
If the eigenvalues are not degenerated,
the eigenvectors are uniquely defined.
Otherwise, we need some additional considerations.

Taking the three principal axes as the axes for a new coordinate system,
the inertia tensor has 
the diagonalized
form
\begin{equation}\label{Inertia-diag}
I=\left(\begin{array}{ccc}
\sum_{k} m_k (Y_k^2+Z_k^2) & 0 & 0 \\
0 & \sum_{k} m_k (Z_k^2+X_k^2) & 0 \\
0 & 0 & \sum_{k} m_k (X_k^2+Y_k^2)
\end{array}\right).
\end{equation}

The condition $I_{xz}=I_{yz}=0$
is identical to the equation $c_x=c_y=0$ for $\omega\ne 0$
in equation (\ref{CxCy}).
Indeed, if (\ref{CxCy}) is satisfied, then the inertia tensor
is partially diagonalized in the form. 

\begin{equation}\label{Inertia-par}
I=\left(\begin{array}{ccc}
I_{xx} & I_{xy} & 0\\
I_{yx} & I_{yy} & 0 \\
0 & 0 & I_{zz}
\end{array}\right).
\end{equation}
The meaning of equation (\ref{CxCy}) is now clear.
That is, the $z$-axis is one of the principal axes of the inertia tensor.
Therefore, the problem of finding a candidate for the z-axis through a given shape variable $\sigma_{ij}$ comes down to computing the inertia tensor and solving its eigenvalue problem.

If a shape has symmetry, then
it is simple to calculate the inertia tensor
and diagonalize it to find the $z$--axis. The simple examples are
the Eulerian relative equilibria.
We will show this in the next subsection.

\subsection{Classification for the relative equilibria}
\begin{figure}
   \centering
   \includegraphics[width=5cm]{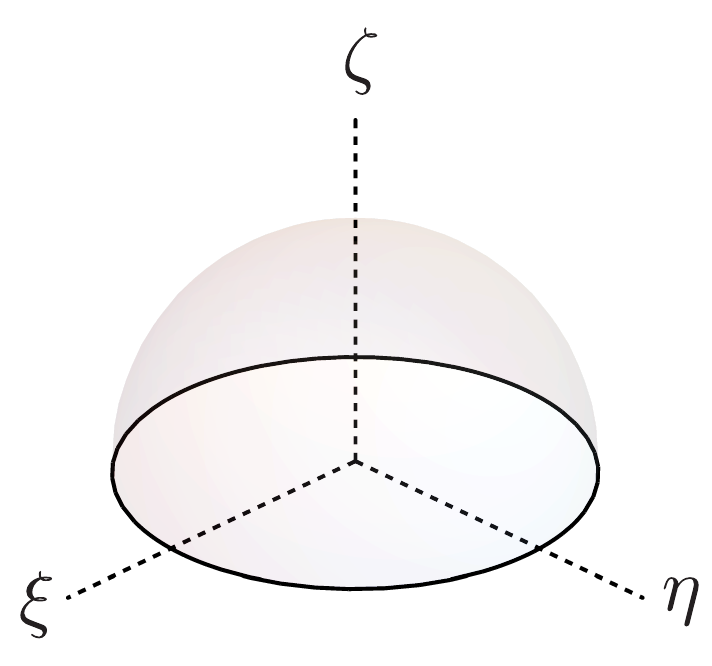} 
   \caption{The masses are on one geodesic
   (black cirlce),
   located on the $\zeta=0$ plane.}
   \label{figPutTheRing}
\end{figure}
For the Eulerian relative equilibria,
we can put the geodesic on the $\zeta=0$ plane (see Figure \ref{figPutTheRing}).
Then, $I_{\xi\zeta}=I_{\eta\zeta}=0$.
To avoid possible confusions
with coordinate $XYZ$,
we use a new $\xi\eta\zeta$ coordinate system for temporal use.

The obvious choice for the principal axis is the $\zeta$--axis.
Taking $\zeta$--axis for $z$--axis,
the geodesic is the Equator.
This type of Eulerian relative equilibria
are  relative equilibria on the Equator.
For this case, $\theta_k=\pi/2$ for $k=1,2,3$ and
to find the relative equilibrium is
reduced to find $\phi_i-\phi_j$ (see \cite{pre} for details).

The other principal axes are on the $\xi\eta$ plane,
because they must be orthogonal to the $\zeta$--axis.
Taking the $z$--axis on the $\xi\eta$ plane,
the geodesic passes through the $z$--axis,
therefore through the north and the south 
pole of the sphere $\mathbb{S}^2$.
So these type of Eulerian relative equilibria
are  relative equilibria on a rotating meridian \cite{pre}. 

Since the extended Lagrangian relative equilibria
are not on a geodesic,
we have to solve the eigenvalue problem for them.
This will be shown in the following subsection. 

\subsection{The inertia tensor
for the extended Lagrangian relative equilibria}
A triangle
with given masses $m_k$ and
given arc angles $\sigma_{ij}$,
can be put in the temporal position
with 
$(\theta_3,\phi_3)=(0,0)$,
$(\theta_1,\phi_1)=(\sigma_{31},0)$,
and $(\theta_2,\phi_2)=(\sigma_{23},\alpha)$ (See Figure \ref{figPutTheTriangle})
where
\begin{equation}
\cos\alpha
=\frac{\cos\sigma_{12}-\cos\sigma_{31}\cos\sigma_{23}}
	{\sin\sigma_{31}\sin\sigma_{23}}.
\end{equation}

\begin{figure}
   \centering
   \includegraphics[width=5cm]{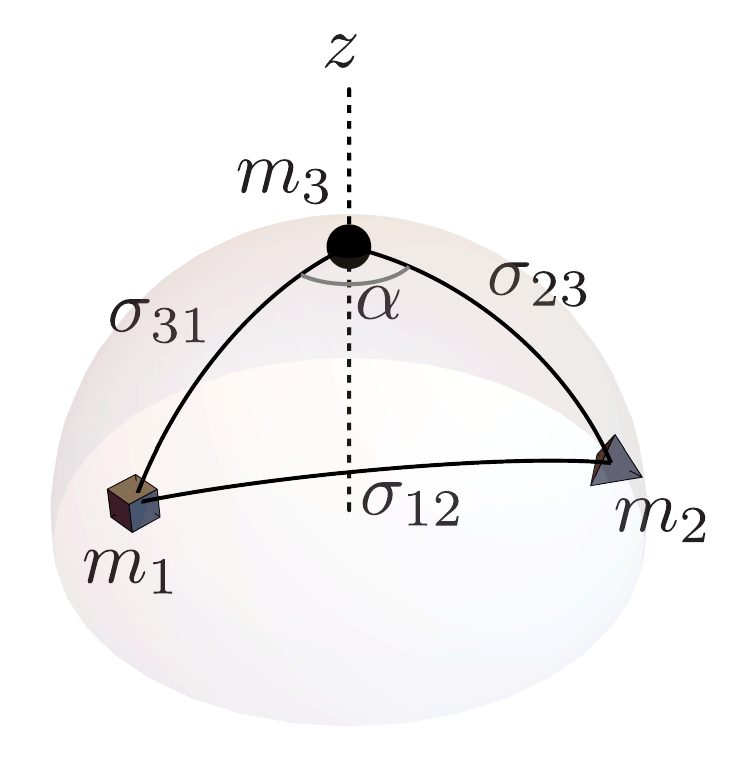} 
   \caption{The triangle with
   $\sigma_{ij}$ is placed to the temporal position
   $(\theta_3,\phi_3)=(0,0)$,
   $(\theta_1,\phi_1)=(\sigma_{31},0)$,
   and
    $(\theta_2,\phi_2)=(\sigma_{23},\alpha)$.
    The cube, tetrahedron, and ball represent
    $m_1$, $m_2$, and $m_3$ respectively. We will continue using this convention in the next figures.}
   \label{figPutTheTriangle}
\end{figure}

In this setting, the components of the inertia tensor are
\begin{equation}
\label{Idiagonal2}
\begin{split}
I_{xx}
&=m_1\cos^2(\sigma_{31})
	+m_2\Big(\cos^2(\sigma_{23})+\sin^2(\sigma_{23})\sin^2(\alpha)\Big)
	+m_3,\\
I_{yy}
&=m_3+m_1+m_2\Big(\cos^2(\sigma_{23})+\sin^2(\sigma_{23})\cos^2(\alpha)\Big),\\
I_{zz}
&=m_1\sin^2(\sigma_{31})+m_2\sin^2(\sigma_{23}),
\end{split}
\end{equation}
and
\begin{equation}
\label{Ioffdiagonal2}
\begin{split}
I_{xy}&=I_{yx}=-m_2\sin^2(\sigma_{23})\sin\alpha\cos\alpha,\\
I_{yz}&=I_{zy}=-m_2\sin\sigma_{23}\cos\sigma_{23}\sin\alpha,\\
I_{zx}&=I_{xz}=-m_1\sin\sigma_{31}\cos\sigma_{31}
		-m_2\sin\sigma_{23}\cos\sigma_{23}\cos\alpha.
\end{split}
\end{equation}

Although the expressions for $I$ are not symmetric for the indexes of masses, the characteristic polynomial is symmetric,
and depends only on the masses $m_k$ and 
the arc angles $\sigma_{ij}$, as it should be.
It has the form
\begin{equation}
\label{defILambda}
\begin{split}
p_{I}(\lambda)=&|\lambda -I|\\
=&\Big(\lambda-(m_1+m_2)\Big)
	\Big(\lambda-(m_2+m_3)\Big)
	\Big(\lambda-(m_3+m_1)\Big)\\
	&-\Big(\lambda-(m_1+m_2)\Big)m_1m_2\cos^2(\sigma_{12})\\
	&-\Big(\lambda-(m_2+m_3)\Big)m_2m_3\cos^2(\sigma_{23})\\
	&-\Big(\lambda-(m_3+m_1)\Big)m_3m_1\cos^2(\sigma_{31})\\
&+2m_1m_2m_3\cos\sigma_{12}\cos\sigma_{23}\cos\sigma_{31}.
\end{split}
\end{equation}

We can calculate the eigenvalue problem \eqref{eigenvalue-problem},
and find the three principal axes.
Choosing one of the principal axis as the $z$--axis,
we can calculate $\cos\theta_k$
which is the angle of the position of $m_k$ from the $z$--axis.

This could be a tedious problem, but fortunately, we have found an alternative expression $J$ for the inertia tensor. 
This is symmetric with respect to the mass index and gives us the $z$--axis and $\cos\theta_k$ simultaneously.

To get $J$, let us start from the equations
(\ref{fundamentalrelation}) and (\ref{CxCy}).
We calculate
\begin{equation}
\begin{split}
&(m_i+m_j)\cos\theta_k
-\cos(\sigma_{ki})m_i \cos\theta_i
-\cos(\sigma_{kj})m_j \cos\theta_j\\
=&(m_i+m_j)\cos\theta_k\\
	&-\big(\cos\theta_k\cos\theta_i+\sin\theta_k\sin\theta_i\cos(\theta_k-\theta_i)\big)m_i \cos\theta_i\\
	&-\big(\cos\theta_k\cos\theta_j+\sin\theta_k\sin\theta_j\cos(\theta_k-\theta_j)\big)m_j \cos\theta_j\\
=&\cos\theta_k
	\big(m_i+m_j-m_i\cos^2(\theta_i)-m_j\cos^2(\theta_j)\big)\\
	&-\sin\theta_k
	\big(
		m_i\sin\theta_i\cos\theta_j\cos(\phi_k-\phi_i)
		+m_j\sin\theta_j\cos\theta_j\cos(\phi_k-\phi_j)
	\big)\\
=&\cos\theta_k
	\big(m_i\sin^2(\theta_i)+m_j\sin^2(\theta_j)\big)\\
	&+\sin\theta_k \big(m_k\sin\theta_k\cos\theta_k\big)
	\mbox{ (here we have used (\ref{cxcyInE3}))}\\
=&\cos\theta_k \sum_{l=1,2,3} m_l\sin^2(\theta_l).
\end{split}
\end{equation}
Multiplying by $\sqrt{m_k}$, 
both the first line and the last line of the above 
equality, we obtain
\begin{equation}
\begin{split}
&(m_i+m_j)\,\sqrt{m_k}\cos\theta_k
-\sqrt{m_km_i}\cos\sigma_{ki}\,\sqrt{m_i}\cos\theta_i
-\sqrt{m_km_j}\cos\sigma_{kj}\,\sqrt{m_j}\cos\theta_j\\
&=\,\sqrt{m_k}\cos\theta_k
	\sum_{l=1,2,3} m_l\sin^2(\theta_l).
\end{split}
\end{equation}
Assembling these equations in a matrix form,
we obtain
\begin{equation}
\begin{split}
&\left(\begin{array}{ccc}
m_2+m_3 & -\sqrt{m_1m_2}\cos\sigma_{12} & -\sqrt{m_1m_3}\cos\sigma_{13} \\
-\sqrt{m_2m_1}\cos\sigma_{21} &m_3+m_1& -\sqrt{m_2m_3}\cos\sigma_{23} \\
-\sqrt{m_3m_1}\cos\sigma_{31} & -\sqrt{m_3m_2}\cos\sigma_{32} & m_1+m_2
\end{array}\right)
\left(\begin{array}{c}
\sqrt{m_1}\cos\theta_1 \\
\sqrt{m_2}\cos\theta_2 \\
\sqrt{m_3}\cos\theta_3
\end{array}
\right)\\
&=\left(\sum_{l=1,2,3} m_l\sin^2(\theta_l)\right)
\left(\begin{array}{c}
\sqrt{m_1}\cos\theta_1 \\
\sqrt{m_2}\cos\theta_2 \\
\sqrt{m_3}\cos\theta_3
\end{array}
\right).
\end{split}
\end{equation}
This means that
the vector $^t(\sqrt{m_1}\cos\theta_1,
\sqrt{m_2}\cos\theta_2,
\sqrt{m_3}\cos\theta_3)$
is an eigenvector of the symmetric matrix
\begin{equation}
\label{defJ}
J=\left(\begin{array}{ccc}
m_2+m_3 & -\sqrt{m_1m_2}\cos\sigma_{12} & -\sqrt{m_1m_3}\cos\sigma_{13} \\
-\sqrt{m_2m_1}\cos\sigma_{21} &m_3+m_1& -\sqrt{m_2m_3}\cos\sigma_{23} \\
-\sqrt{m_3m_1}\cos\sigma_{31} & -\sqrt{m_3m_2}\cos\sigma_{32} & m_1+m_2
\end{array}\right),
\end{equation}
with the eigenvalue $\lambda=I_{zz}$.
This eigenvalue problem  gives directly $\cos\theta_k$.
To give a correct value of $\cos\theta_k$,
the eigenvector must be normalised to be
\begin{equation}
\sum_k (\sqrt{m_k}\cos\theta_k)^2
=M-\sum_k m_k \sin^2(\theta_k)
=M-\lambda,
\end{equation}
where $M=\sum_k m_k$ is the total mass. 

Therefore, for the eigenvalue $\lambda$ and 
the normalised eigenvector $\Psi=(\psi_1,\psi_2,\psi_3)$,
the angle $\theta_{k}$ is given by
\begin{equation}
\label{cosTheta}
\cos\theta_k=\frac{\psi_k \sqrt{M-\lambda}}{\sqrt{m_k}} \quad
\mbox{ with } \quad
\sum_k \psi_k^2=1.
\end{equation}

The characteristic polynomial for $J$
is exactly the same as that for the inertia tensor $I$
in (\ref{defILambda}).
Since both $I$ and $J$ are representations of symmetric matrices,
they must be similar matrices,
\begin{equation}
^tQIQ=J,
\quad
^tQ=Q^{-1},
\end{equation}
where the matrix $Q$ is orthogonal.
Namely, the difference between $I$ and $J$
is the placement of the same shape $\sigma_{ij}$.
In other words,
the matrix $J$ is another expression for the inertia tensor $I$.
Therefore,
the eigenvectors for $J$ give  the principal axes for 
the shape $\sigma_{ij}$.

Since $J$ depends only on $m_k$ and $\sigma_{ij}$,
the quantities $\cos\theta_k$ are given as a function of them.
Thus we obtain 
the translation formula for
a given shape
$\sigma_{ij}$ to the angle $\theta_k$.
Utilising this correspondence,
the condition (\ref{conditionForLagrangian4})
can be written as
\begin{equation}
-\frac{2}{\omega^2}\sum_k m_k \cos^2(\theta_k)
=\frac{\cos\theta_i\cos\theta_j}{U'(D_{ij}^2)}
=\frac{\psi_i\psi_j}{\sqrt{m_im_j}}\frac{(M-\lambda)}{U'(D_{ij}^2)}.
\end{equation}
Since, the eigenvalue $\lambda=\sum_k m_k \sin^2(\theta_k)$,
we obtain that $M-\lambda=\sum_k m_k \cos^2(\theta_k)$
and it is not zero for the Lagrangian case, we finally obtain
the condition for a rigid rotator
\begin{equation}
\label{conditionForRigidRotator0}
-\frac{2}{\omega^2}
=\frac{\psi_i\psi_j}{\sqrt{m_im_j}\,\,U'(D_{ij}^2)}.
\end{equation}

For the particular case of the cotangent potential the condition is
\begin{equation}
\label{conditionForRigidRotator}
\frac{1}{R^3\omega^2}
=\frac{\psi_i\psi_j}{\sqrt{m_im_j}}\sin^3(\sigma_{ij})
\end{equation}
for all $i,j$.

In this way we have gotten a method to find the Lagrangian
rigid rotators.
Then, by the equation (\ref{cosTheta}),
the quantities $\cos\theta_k$ are determined.
With this we get the Lagrangian $RE$ configuration.
We will show some interesting rigid rotators
in the following section.


\section{Some extended Lagrangian relative equilibria}\label{LRE}
In 2012 Diacu, Per\'{e}z-Chavela, and Santoprete proved that
an equilateral triangular rigid  rotator is feasible 
only when the masses are equal \cite{Diacu-EPC1}.
However, a rigid  rotator of equal masses 
does not necessarily have to be an equilateral triangle.
In fact, there are isosceles triangle rigid rotators other than equilateral triangles.

In the following subsections, 
first we will prove the result of Diacu et. al. in our context,
actually we have generalized this result for generic potentials.
In the other subsections, in order to have concrete results we have restricted our analysis to the cotangent potential.
These results give us new families of Lagrangian $RE$ for the two dimensional positive curved three body problem. As far as we know, this is the first time such families are shown.

First we study the case of   
an almost trivial isosceles triangle with equal masses,
then general isosceles triangles with equal masses and with only two equal masses.

\subsection{Equilateral extended Lagrangian rigid rotator
for generic potential}
In this subsection,
we show that
an equilateral triangle can be an extended
Lagrangian rigid rotator
if and only if the masses are equal and the angles $\theta_k$ are equal.
This statement is true for generic potential $U(D^2)$.

For $\sigma_{12}=\sigma_{23}=\sigma_{31}$,
$D_{12}^2=D_{23}^2=D_{23}^2$,
therefore, $U'(D_{12}^2)=U'(D_{23}^2)=U'(D_{31}^2)$.
Then, by the condition (\ref{conditionForLagrangian2}),
$\cos\theta_1=\cos\theta_2=\cos\theta_3$.
Therefore, the matrix $J$ must have the eigenvector 
$^t\!(\sqrt{m_1},\sqrt{m_2},\sqrt{m_3})$. 
\begin{equation}
\left(\begin{array}{ccc}
m_2+m_3 & -\sqrt{m_1m_2}\cos\sigma & -\sqrt{m_1m_3}\cos\sigma \\
-\sqrt{m_2m_1}\cos\sigma &m_3+m_1& -\sqrt{m_2m_3}\cos\sigma \\
-\sqrt{m_3m_1}\cos\sigma & -\sqrt{m_3m_2}\cos\sigma & m_1+m_2
\end{array}\right)
\left(\begin{array}{c}
\sqrt{m_1}\\\sqrt{m_2}\\\sqrt{m_3}
\end{array}\right)
=\lambda
\left(\begin{array}{c}
\sqrt{m_1}\\\sqrt{m_2}\\\sqrt{m_3}\end{array}\right).
\end{equation}
Where $\sigma=\sigma_{12}=\sigma_{23}=\sigma_{31}$.
This equation reduces to
\begin{equation}
\lambda
=(m_1+m_2)(1 - \cos\sigma)
=(m_2+m_3)(1 - \cos\sigma)
=(m_3+m_1)(1 - \cos\sigma),
\end{equation}
which means that $m_1=m_2=m_3$ because $\cos\sigma\ne 1$. 

\subsection{%
Equal masses isosceles triangle with
$\sigma_{23}=\sigma_{31}=\pi/2$  for the cotangent potential}

For  $m_1=m_2=m_3=m$, $\sigma_{23}=\sigma_{31}=\pi/2$,
and $\sigma_{12}=\sigma$,  
the matrix $J$ takes the form
\begin{equation}
J
=\left(\begin{array}{ccc}
2m & -m\cos\sigma & 0 \\
-m\cos\sigma & 2m & 0 \\
0 & 0 & 2m\end{array}
\right).
\end{equation}
The eigenvalues and eigenvectors are
\begin{align}
&\lambda_0=2m,&{}^t(0,0,1),\\
&\lambda_\pm=m(2\pm\cos\sigma),&{}^t(\mp1,1,0).
\end{align}

Therefore, if $\sigma\ne \pi/2$,
the eigenvalues and eigenvectors are different each other.
Obviously, there are no eigenvectors for which all components have the same sign.
By Proposition \ref{cosMustHaveTheSameSign},
all $\cos\theta_k$ must have the same sign, so all components of the eigenvectors
must do.

When $\sigma=\pi/2$,
the eigenvalues are triply degenerate,
for this case,
we can take the eigenvector ${}^t(1,1,1)$.
This corresponds to the equilateral solution $\sigma_{ij}=\pi/2$ for all $i,j$.
In this case, $\lambda=2m$, the normalised eigenvector is 
$\Psi=(\psi_1,\psi_2,\psi_3)=\nobreak{^t(1,1,1)/\sqrt{3}}$.
Obviously, the condition for a rigid rotator is satisfied,
\begin{equation}
\frac{1}{R^3\omega^2}=\frac{\psi_i\psi_j}{m}\sin^3(\sigma_{ij})=\frac{1}{3m}.
\end{equation}
Therefore, this equilateral triangle is a rigid rotator.
The values for this solution are
$\cos\theta_k=\psi_k\sqrt{(M-\lambda)/m}=1/\sqrt{3}$,
$\sin\theta_k=\sqrt{2/3}$, and $\cos(\phi_i-\phi_j)=-1/2$.

\subsection{%
Equal masses isosceles triangle with $\sigma_{23}=\sigma_{31}$ for the cotangent potential}
\label{secEqualMassIsosceles}

Now, we consider general equal masses isosceles triangle,
$m_1=m_2=m_3=m$, and
$\sigma_{23}=\sigma_{31}=\sigma$.

In this case the matrix $J$ is given by
\begin{equation}
J=
\left(\begin{array}{ccc}
2m & -m\cos\sigma_{12} & -m\cos\sigma \\
-m\cos\sigma_{12} & 2m &  -m\cos\sigma \\
 -m\cos\sigma &  -m\cos\sigma & 2m
\end{array}\right).
\end{equation}

The eigenvalues and the 
un-normalised eigenvectors are
\begin{equation}
\lambda_0=(2+\cos\sigma_{12})m,\,
\psi_0
={}^t(-1,1,0),
\end{equation}
\begin{equation}
\label{eigenFunctionsForEqualMassIsosceles}
\begin{split}
\lambda_\pm
&=\frac{m}{4}\left(
	8-2\cos\sigma_{12}
	\pm 2\sqrt{8\cos^2(\sigma)+\cos^2(\sigma_{12})}
	\right),\\
\psi_\pm
&=
\left(\begin{array}{c}
\left(2\cos\sigma_{12}\mp
	2\sqrt{8\cos^2(\sigma)+\cos^2(\sigma_{12})
	}\right)/\cos\sigma\\
\left(2\cos\sigma_{12}\mp
	2\sqrt{8\cos^2(\sigma)+\cos^2(\sigma_{12})}
	\right)/\cos\sigma \\
8
\end{array}\right),
\end{split}
\end{equation}
for $\cos\sigma\ne 0$.
The eigenfunctions $\psi_\pm$ are not  properly defined
for $\cos\sigma=0$. In the previous subsection,
we know that the solution for $\sigma=\pi/2$ is just
$\sigma_{12}=\pi/2$.
So, in this subsection, we exclude $\sigma=\pi/2$.

The eigenvalues are triply degenerate
if and only if $\sigma=\sigma_{12}=\pi/2$.
Since
we are considering the case $\sigma\ne \pi/2$,
the three eigenvalues are different from each other.

For  each normalised components $\psi_i$ of
the eigenfunction $\psi_\pm$,
the equation $\psi_i\psi_j\sin^3(\sigma_{ij})=$ common, with $\sigma\ne \pi/2$ 
reduces to
\begin{equation}
\label{conditionForEqualMassEquilateral}
\begin{split}
-4\sin^3(\sigma)\cos\sigma+\sin^3(\sigma_{12})\cos\sigma_{12}
=\pm \sin^3(\sigma_{12})
\sqrt{8\cos^2(\sigma)+\cos^2(\sigma_{12})}.\\
\end{split}
\end{equation}
If we square the above equation, we obtain
\begin{equation}
32\cos\sigma
\Big(
\cos\sigma
\big(2\sin^6(\sigma)-\sin^6(\sigma_{12})\big)
-\sin^3(\sigma)\cos(\sigma_{12})\sin^3(\sigma_{12})
\Big)=0.
\end{equation}
Since, we are considering $\cos\sigma\ne 0$,
we get
\begin{equation}
q(\sigma,\sigma_{12})
=\cos\sigma
\big(2\sin^6(\sigma)-\sin^6(\sigma_{12})\big)
-\sin^3(\sigma)\cos(\sigma_{12})\sin^3(\sigma_{12})
=0.
\end{equation}

The graphical representation of this equation is
shown in Figure \ref{figMspecialCase8equalMassesCaseCont3FigSolution}.
\begin{figure}
   \centering
   \includegraphics[width=7cm]{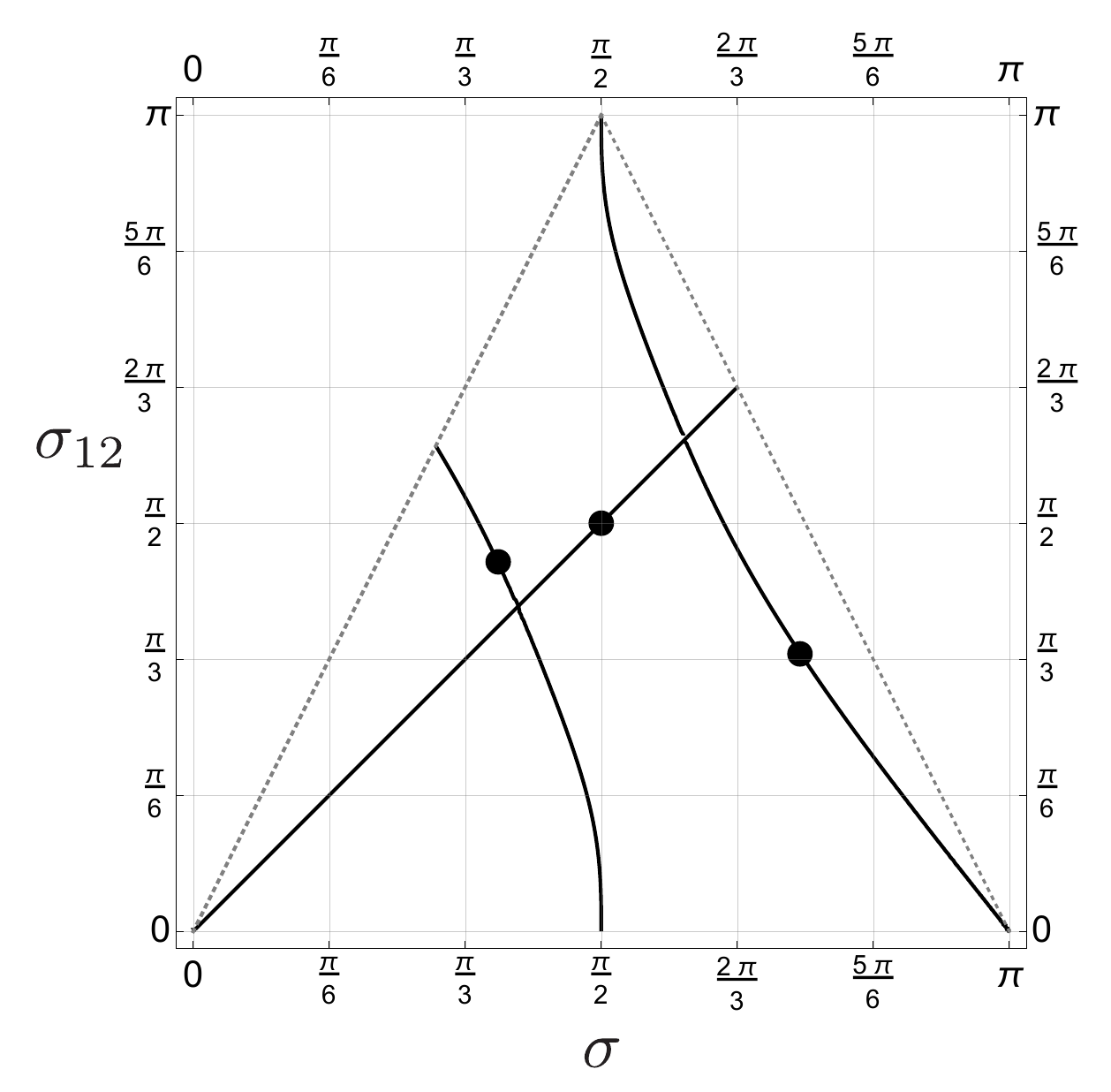} 
   \caption{The point $(\sigma,\sigma_{12})$ on the  solid line represents the corresponding equal masses isosceles 
   rigid rotator
   with $\sigma_{12}$ and $\sigma=\sigma_{23}=\sigma_{31}$.
   The straight line represents
   equilateral rigid rotators.
   The region inside the grey lines represents 
   the region to form a triangle, 
   $\sigma_{12}<2\sigma<2\pi-\sigma_{12}$.
   Note that the curve is point symmetric
   around $(\sigma,\sigma_{12})=(\pi/2,\pi/2)$.
   The three  black circles represent
   the right-angled rigid rotators,
   where the angle at the vertex
   $m_3$ is $\pi/2$.
   }
   \label{figMspecialCase8equalMassesCaseCont3FigSolution}
\end{figure}

The most important aspect on the figure 
is that the curves
are point symmetric around $(\sigma,\sigma_{12})=(\pi/2,\pi/2)$.
Indeed, for the transformation
\begin{equation}
\label{theTransformation}
\sigma\to \pi-\sigma,\,\,
\sigma_{12}\to \pi-\sigma_{12},
\end{equation}
the left hand side of (\ref{conditionForEqualMassEquilateral})
is invariant,
while
the right hand side jumps to the opposite sign branch,
namely, $+$ branch to $-$ and vice versa.
This is also true for the eigenvectors $\psi_\pm$ in (\ref{eigenFunctionsForEqualMassIsosceles}).
Therefore, the transformation (\ref{theTransformation})
maps a solution to another solution
if the arcs can form a triangle,
in other words, if the mapped point is inside the grey lines.
This transformation keeps the angular velocity 
$R^3\omega^2$ invariant. (See Figure \ref{figisoscelesSigma12EqPiDiv3and2PiDiv3}.)
\begin{figure}
   \centering
\includegraphics[width=4cm]{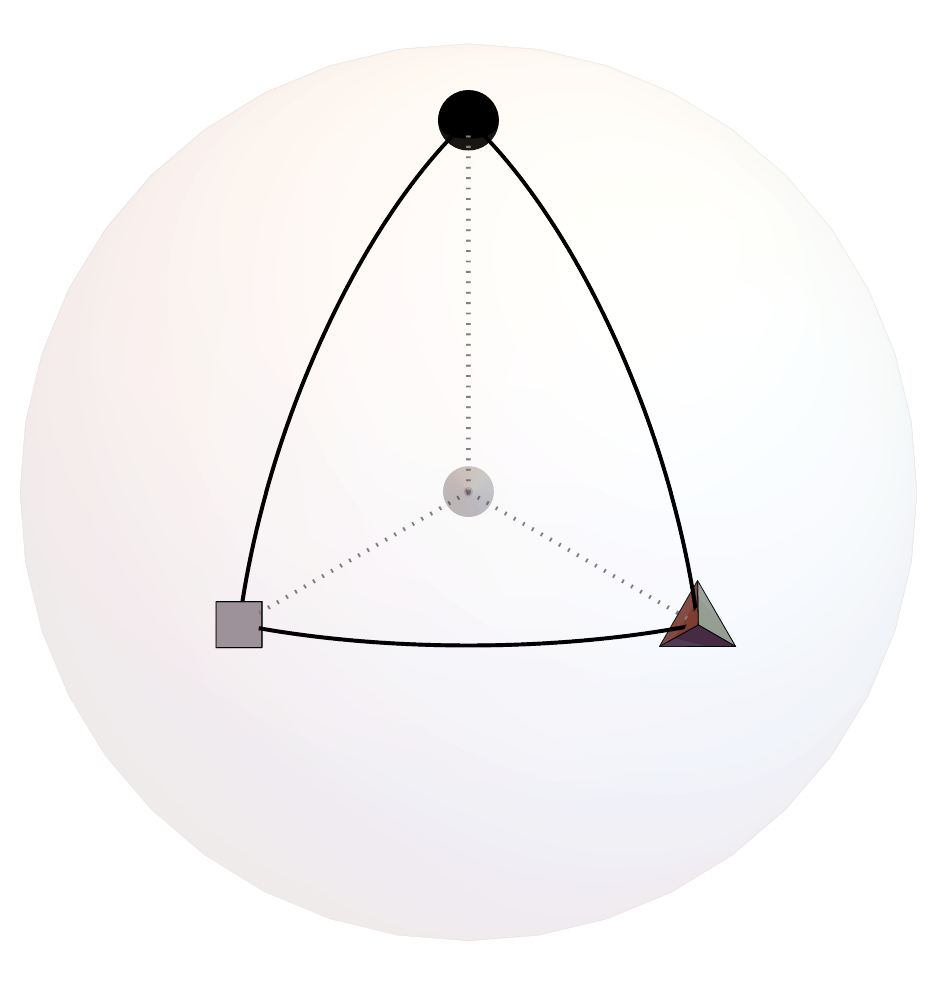}%
\includegraphics[width=4cm]{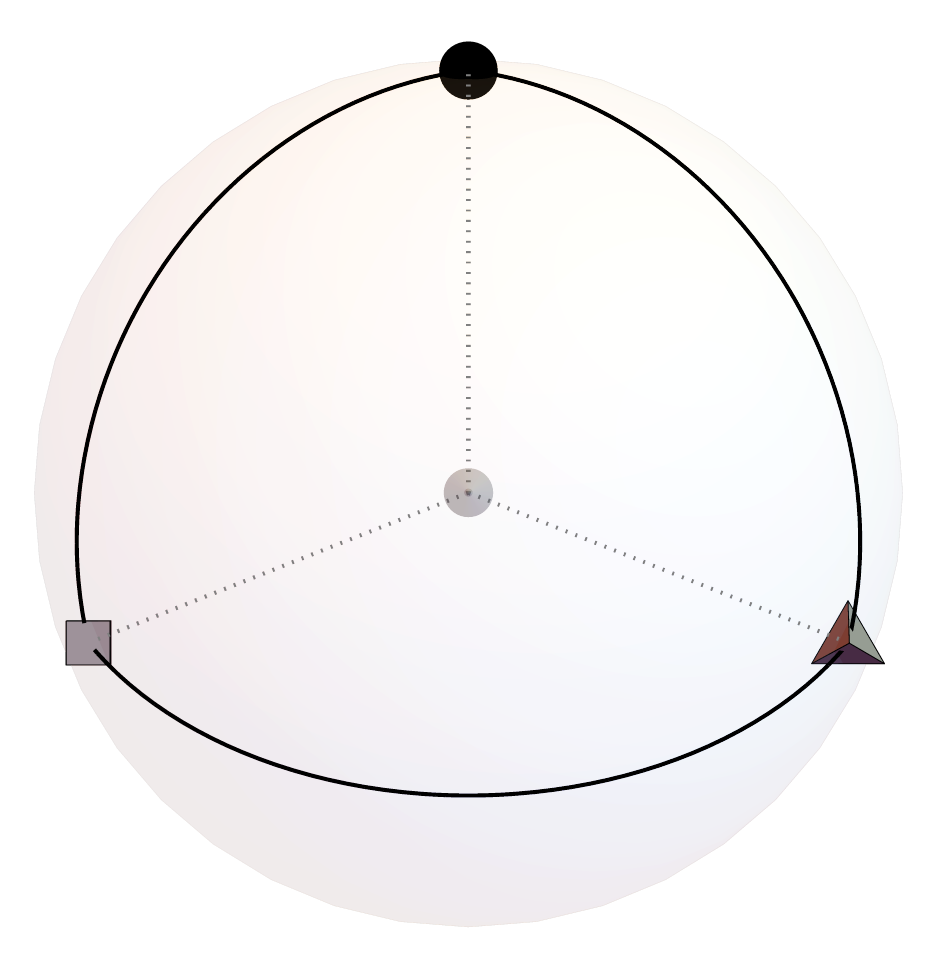}
   \caption{
   Two isosceles Lagrangian $RE$
   configurations
   with $\sigma_{12}=\pi/3$ (left) and $2\pi/3$ (right)
   for equal masses,
   seen from above the North pole.
   The grey ball at the centre represents the North pole.
   Each triangle has
   $\sigma_{23}=\sigma_{31}=1.33240...$ (left) and $1.80918...$ (right)
   respectively.
  Sum of the corresponding arc angle $\sigma_{ij}$ is $\pi$.
   They have the common angular velocity $R^3\omega^2=3.85072...$.}
   \label{figisoscelesSigma12EqPiDiv3and2PiDiv3}
\end{figure}

Other important aspects of Figure \ref{figMspecialCase8equalMassesCaseCont3FigSolution}
are the following: 
The straight line $\sigma=\sigma_{12}$ corresponds to 
equilateral triangles.
The right end of this line is $\sigma=\sigma_{12}=2\pi/3$,
where the triangle becomes  the rigid rotator on 
the Equator.

The crossing points are the saddle points of $q(\sigma,\sigma_{12})$.
These points are the solutions of
\begin{equation}
\frac{\partial q}{\partial \sigma}
=\frac{\partial q}{\partial \sigma_{12}}=0
\mbox{ and } \sigma=\sigma_{12}.
\end{equation}
The solutions are
$\cos\sigma=\cos\sigma_{12}=\pm \sqrt{1/10}$.
Let the smaller value be 
\begin{equation}
\sigma_s=\arccos(\sqrt{1/10})=1.24904...
\end{equation}
Then the larger one is $\pi-\sigma_s=1.89254...$.

The left end of the left curve is the solution of
$q(\sigma,2\sigma)=0$, which is the solution of
\begin{equation}
32\cos^6(\sigma)+8\cos^4(\sigma)-4\cos^2(\sigma)-1=0.
\end{equation}
Let the solution be $\sigma_E$,
\begin{equation}
\begin{split}
&\sigma=\sigma_E=\arccos(2^{-3/4})=0.93402....,\\
&\sigma_{12}=2\sigma_E=1.86804....
\end{split}
\end{equation}
This end point is also an Eulerian rigid rotator.

The Figure \ref{figMspecialCase8equalMassesCaseCont3FigSolution}
also shows that
three solutions exist for $0<\sigma_{12}<2\sigma_E$
with the exceptional point
$\sigma_{12}=\sigma=\sigma_s$.
(See Figure \ref{figIsoscelesSigma12EqPiDiv6}.)
Then for $2\sigma_E<\sigma_{12}<2\pi/3$,
two solutions exist
with exceptional point $\sigma_{12}=\sigma=\pi-\sigma_s$.
Finally, in $2\pi/3<\sigma_{12}<\pi$,
 only one solution exists.
 \begin{figure}
    \centering
\includegraphics[width=4cm]
{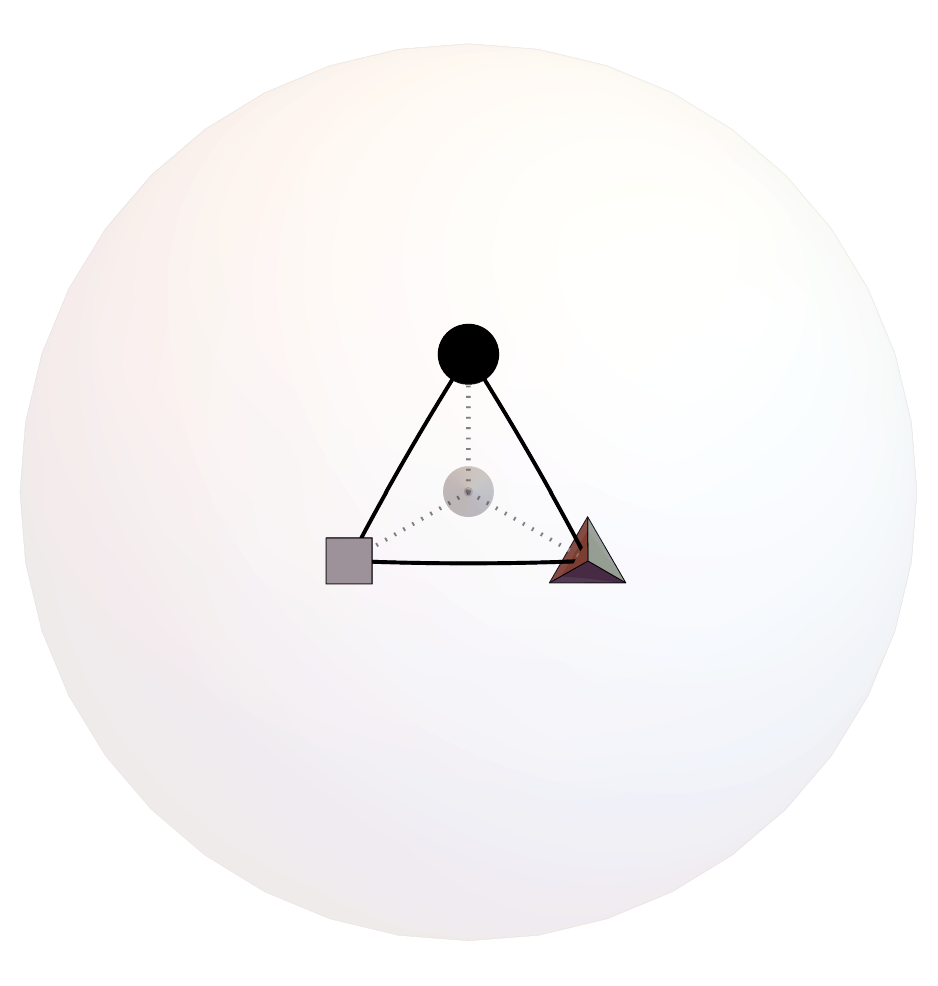}%
\includegraphics[width=4cm]{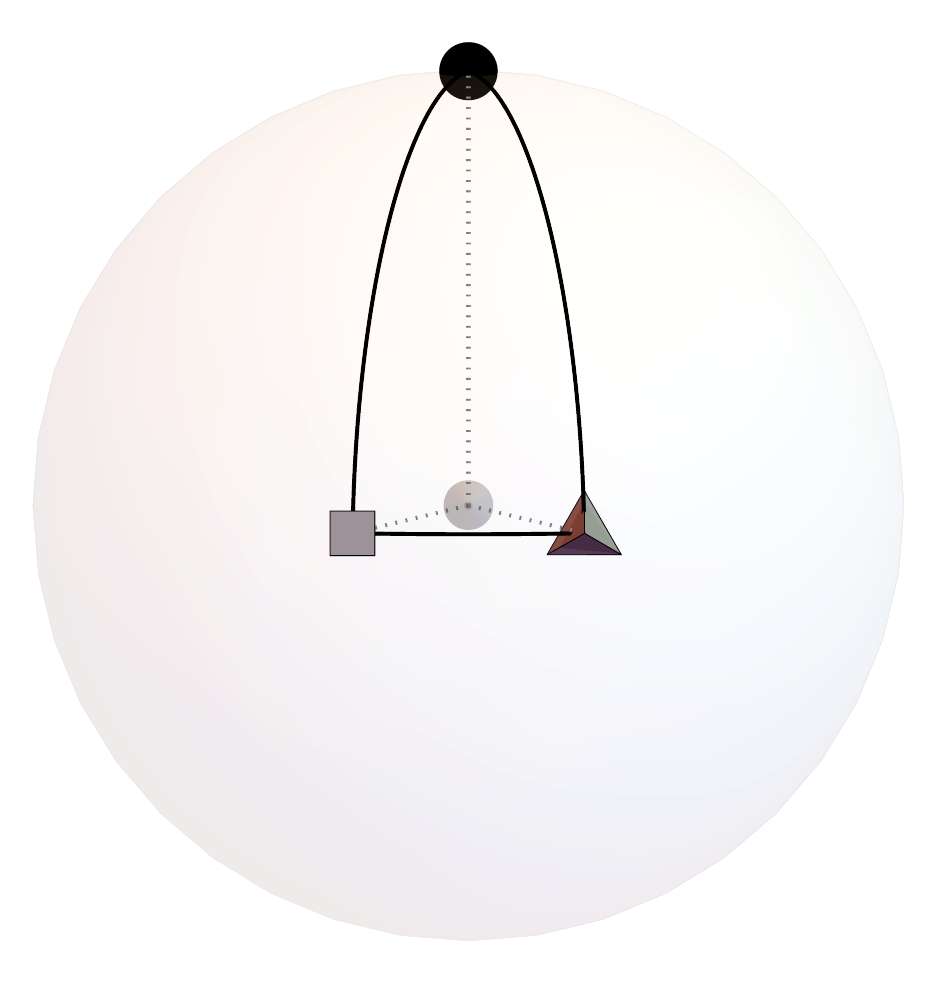}%
\includegraphics[width=4cm]{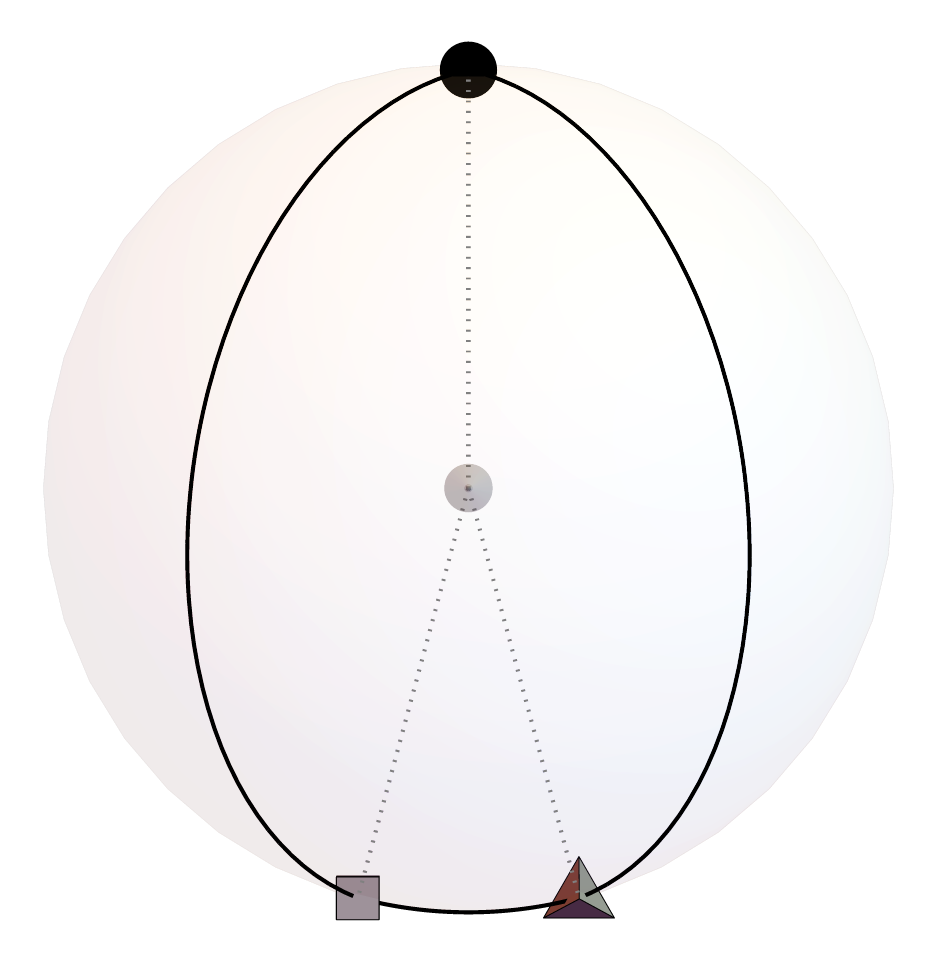} 
    \caption{Three isosceles Lagrangian $RE$ configurations
    with $\sigma_{12}=\pi/6$
    for equal masses,
    seen from above the North pole.
    From left to 
    right
    $\sigma_{23}=\sigma_{31}=
    \pi/6$ (equilateral triangle), $1.51596...$, and $2.73083...$.}
    \label{figIsoscelesSigma12EqPiDiv6}
 \end{figure}

The figure \ref{figMspecialCase8equalMassesCaseCont3FigSolution}
contains three right-angled isosceles triangles,
where the angle at the vertex $m_3$ is $\pi/2$.
They are the solutions of $\cos(\sigma_{12})=\cos^2(\sigma)$.
The obvious one is $\sigma_{12}=\sigma=\pi/2$.
For this shape, the angle of three vertexes are $\pi/2$.
The other two are
$\sigma=1.17275...$ and 
$2.33759...$.

\subsection{Isosceles rigid rotators
with $m_1=m_2$,
$\sigma_{12}=\pi/2$ and  
$\sigma_{23}=\sigma_{31}$ for the cotangent potential}

In the previous subsection,
it was shown that
there are three equal masses rigid rotators
with $\sigma_{12}=\pi/2$.
Therefore, for $\sigma_{12}=\pi/2$, one would expect 
a rigid rotator to exist for approximately equal masses.
What happens when the difference in masses is increased?

To give a partial answer to this question,
isosceles triangles with $m_1=m_2=m\nu$, $m_3=m$,
$\sigma_{12}=\pi/2$,  and
$\sigma_{23}=\sigma_{31}=\sigma$
is considered in this subsection.
For this case,
the matrix $J$ takes the form
\begin{equation}
\left(\begin{array}{ccc}
m(1+\nu) & 0 & -m\sqrt{\nu}\cos\sigma \\
0 & m(1+\nu) &   -m\sqrt{\nu}\cos\sigma \\
 -m\sqrt{\nu}\cos\sigma &   -m\sqrt{\nu}\cos\sigma & 2m\nu
\end{array}\right),
\end{equation}
and the characteristic polynomial is
\begin{equation}
p(\lambda)
=\Big(\lambda-m(1+\nu)\Big)
	\Big(
	(\lambda-m\nu)(\lambda-m(1+2\nu))
		-m^2\nu\cos(2\sigma)
	\Big).
\end{equation}
The eigenvalues and un-normalised eigenvectors are
\begin{equation}
\begin{split}
\lambda_0=m(1+\nu),\qquad
\Psi_0={}^t(-1,1,0),
\end{split}
\end{equation}
\begin{equation}
\begin{split}
\lambda_{\pm}
=&\frac{m}{2}\left(
	(1+3\nu)
	\pm\sqrt{(1-\nu)^2+8\nu\cos^2(\sigma)}
	\right),\\
\Psi_{\pm}
=&
\left(\begin{array}{c}
2\sqrt{\nu}\cos\sigma\\
2\sqrt{\nu}\cos\sigma\\
1-\nu\mp\sqrt{(1-\nu)^2+8\nu\cos^2(\sigma)}
\end{array}\right).
\end{split}
\end{equation}

Note that for $\sigma=\pi/2$,

\begin{equation}
\begin{split}
\lambda_\pm
=\frac{m}{2}\Big(1+3\nu\pm|\nu-1|\Big)
\end{split}
\end{equation}
Therefore, the eigenvalue is
doubly degenerate
for $\nu>1$, $\lambda_0=\lambda_-=m(1+\nu)$,
$\lambda_+=2m\nu$, and
for $\nu<1$, $\lambda_0=\lambda_+=m(1+\nu)$,
$\lambda_-=2m\nu$.
For $\nu=1$, the eigenvalues are triply degenerate
$\lambda_0=\lambda_\pm=2m$.

However,  at $\sigma=\pi/2$, the triangle is equilateral.
We have already investigated the equilateral solution.
So, we can exclude $\sigma=\pi/2$ in this section.
In the following, we assume $\sigma\ne \pi/2$,
where the three eigenvalues are different.

The eigenvector $\Psi_0$ is not suitable,
because it contains a zero component.

For the eigenvector $\Psi_\pm$,
the equation
\begin{equation}
\label{eqLagrangian5}
\frac{\psi_i\psi_j}{\sqrt{m_im_j}}\sin^3(\sigma_{ij})
=\frac{1}{R^3\omega^2}
\end{equation}
yields
\begin{equation}
\label{eqForTheSolution}
2\cos(\sigma)
=\sin^3(\sigma)
	\left(
	1-\nu\mp\sqrt{(1-\nu)^2+8\nu\cos^2(\sigma)}
	\right).
\end{equation}

The solution of $\nu$ is
\begin{equation}
\nu=\frac{\cos(\sigma)-\sin^3(\sigma)}
	{\sin^3(\sigma)(2\sin^3(\sigma)\cos(\sigma)-1)}.
\end{equation}
It is plotted in Figure \ref{isoscelesMassRatio}.
\begin{figure}
   \centering
   \includegraphics[width=7cm]{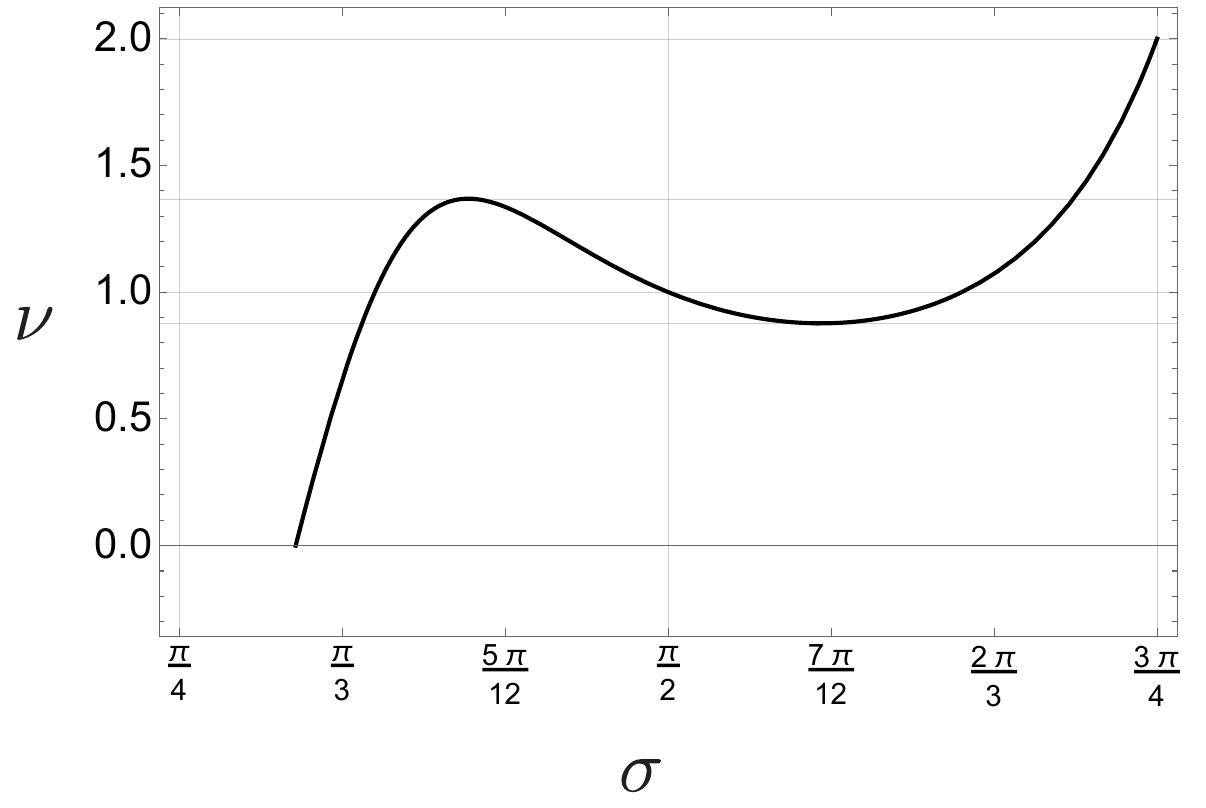} 
   \caption{Mass ratio $\nu=m_1/m_3=m_2/m_3$ for 
   the rigid rotator with 
   $\sigma=\sigma_{23}=\sigma_{31}$, and $\sigma_{12}=\pi/2$.
   }
   \label{isoscelesMassRatio}
\end{figure}

For $0.87... < \nu < 1.36...$, three rigid rotators exist.
For $\nu<0.87...$ or $1.36 < \nu$, only one exists.
The right end point $\nu=2$, $\sigma=3\pi/4$ is
an Eulerian rigid rotator on the Equator.

The left end point is $\nu=0$, 
where $\sigma$ is the solution of $\cos(\sigma_0)=\sin^3(\sigma_0)$.
The numerical value is $\sigma_0=0.97202... > \pi/4$,
The point $\sigma=\pi/4$ is the Eulerian rigid rotator.
Therefore, the left end does not connect to the Eulerian rotator.

Figure \ref{figIsoscelsNuEq1Over100},
shows the Lagrangian $RE$ near the end point,
$\nu=1/100$.
The value of $\sigma$ is $0.97306...>\sigma_0$.

\begin{figure}[htbp] 
   \centering
   \includegraphics[width=4cm]
   {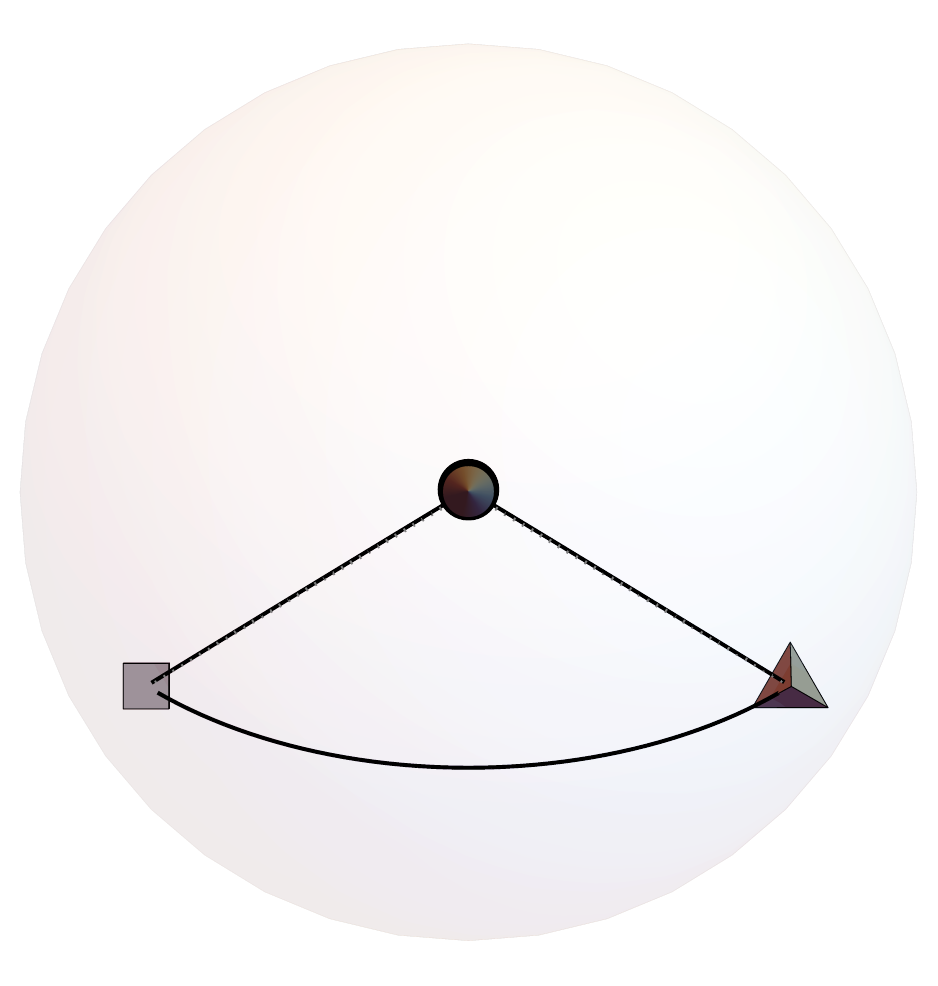} 
   \includegraphics[width=4cm]
   {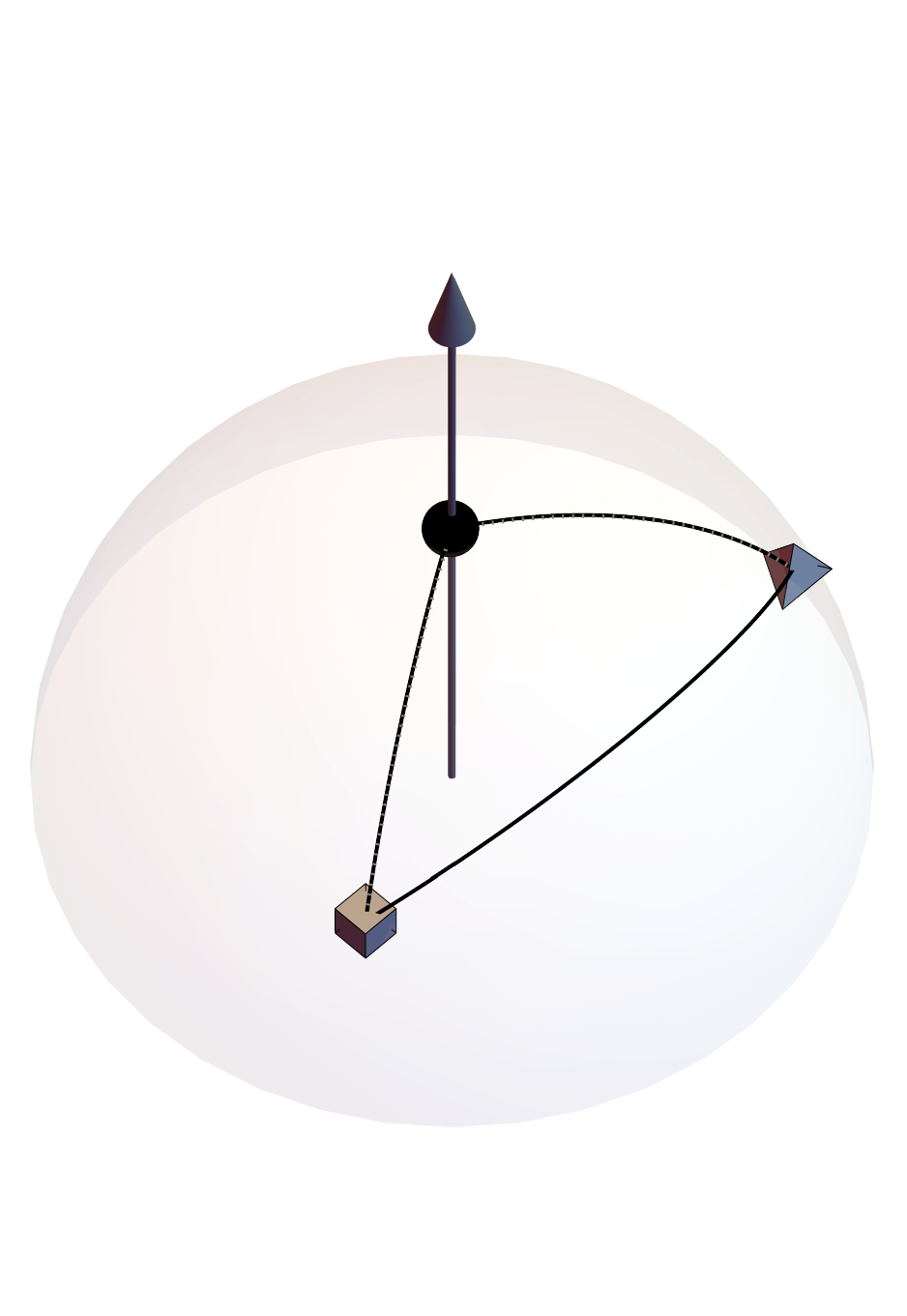} 
  \caption{Isosceles $RE$
  configuration
   with $m_1=m_2=1$, $m_3=100$,
   $\sigma_{12}=\pi/2$,
   $\sigma_{23}=\sigma_{31}=0.97306...$,
   seen from  above the North pole (left) and an arbitrary view point.
   The arrow represents the z-axis.
   The other parameters are
   $R^3\omega^2=315.45...$,
   $\cos\theta_1=\cos\theta_2=0.56481...,
   \cos\theta_3=0.99998...$,
   $\cos(\phi_1-\phi_2)=-0.46846...,
\cos(\phi_2-\phi_3)=\cos(\phi_3-\phi_1)=-0.51552...$.
 $m_3$ is almost on the North pole.
   }
   \label{figIsoscelsNuEq1Over100}
\end{figure}

Although
the limit of the mass ratio $\nu=m_1/m_3=m_2/m_3 \to 0$ is 
a regular point,
$\nu=0$ corresponds to
$\cos\theta_1=\cos\theta_2=\cos\sigma_0$
and $\cos\theta_3=1$,
which does not satisfy the equation of motion.
Actually, $\cos\theta_3=1$ is excluded
for the extended Lagrangian configuration.

\section{Conclusions and final remarks}\label{sec7}
In this paper,
we successfully developed a method to find 
extended Lagrangian relative equilibria on $\mathbb{S}^2$.
The condition for a shape given by $m_k$ and $\sigma_{ij}$
to be a rigid rotator is given by 
the equation 
(\ref{conditionForRigidRotator0}) for a generic potential,
and by equation (\ref{conditionForRigidRotator}) for the cotangent potential. 

When we find a rigid rotator,
we put it on $\mathbb{S}^2$ by using
$\cos\theta_k$ given by the formula (\ref{cosTheta}),
in this way we obtain a Lagrangian $RE$ configuration.

This method for $\mathbb{S}^2$ is parallel to
the method for Euclidian plane. 
The symmetry group of the Euclidian plane is
the three parameter group $T_2\times SO(2)$,
where $T_2$ is the translation
and $SO(2)$ is the rotation.
The $T_2$ invariance yields the conservation of
the centre of mass,
which we take it as the centre for the rotation.
The remaining $SO(2)$ invariance describes
the rotation of the tree bodies around the centre of mass (see the nice Wintner's book for more details \cite{Wintner}).  

Similarly, the symmetry group of $\mathbb{S}^2$ is the
three parameter group $SO(3)$,
that yields the conservation of the angular momentum
${\bf c}=(c_x,c_y,c_z)$.
We take the $z$--axis as the rotation axis by doing $c_x=c_y=0$ .
Then the remaining $SO(2)$ invariance
describes
the rotation of the 
three bodies around the $z$ axis. 

By using this method,
we find some isosceles Lagrangian $RE$ configurations
with three equal masses,
and with just two equal masses.

We observe that even for equal masses, Lagrangian rigid rotators
are non-trivial. It can be equilateral or isosceles triangles.
Then, two questions remain.
\begin{itemize}
\item A scalene triangle (a triangle where the three sides have different lengths) with three equal masses, can be a rigid rotator?
\item  How many rigid rotator exist for given $\sigma_{12}$.
\end{itemize}

The second question is connected with the first one,
because if there are no scalene triangle rigid rotator, then
the maximum number is three as we have shown 
in section \ref{secEqualMassIsosceles}.

In this paper we found just a tiny part of the Lagrangian relative equilibria on $\mathbb{S}^2$. Many equilibria have not yet been found,
and many questions are waiting to be asked.

\appendix

\subsection*{Acknowledgements}
The second author (EPC) has been partially supported 
by Asociaci\'on Mexicana de Cultura A.C. and Conacyt-M\'exico Project A1S10112.


\begin{thebibliography}{99}
\bibitem{Bengochea} Bengochea A., Garc\'ia-Azpeitia C., P\'erez-Chavela E.,  Roldan P. {\it Continuation of relative equilibria in the $n$--body problem to spaces of constant curvature}  
Journal of Differential Equations, \textbf{307}, (2022), 137-159.

\bibitem{Borisov1} Borisov, A.V., Mamaev, I.S., Bisyaev, I.A. {\it Three vortices in spaces of constant curvature: Reduction, Poisson Geometry and Stability}. Regul. Chaot. Dyn. \textbf{23}, (2018), 613-636.


\bibitem{Diacu-EPC1} Diacu F., P\'erez-Chavela E., Santoprete M., {\it The n-body problem in spaces of constant curvature. Part I: Relative equilibria.} J. Nonlinear Sci. \textbf{22} (2012), no. 2, 247--266.

\bibitem{Diacu1} Diacu F., Relative equilibria of the curved N-body problem. Atlantis Studies in Dynamical Systems, Atlantis Press, Amsterdan, Paris, Beijing \textbf{1}, 2012.

\bibitem{Diacu4} Diacu F., S\'anchez-Cerritos J.M. and Zhu S.
\emph{Stability of Fixed Points and Associated Relative Equilibria of the 3-body Problem on $S^1$ and $S^2$}, Journal of Dynamics and Differential Equations, {\bf 30}, (2018), 209-225.

\bibitem{pre} Fujiwara T. and P\'erez-Chavela E. {\it Three body relative equilibria on $\mathbb{S}^2$ I: Euler configurations}. Preprint 2022. 

\bibitem{Goldstein}
H.~Goldstein, C.~Poole, and J.~Safko,
``Classical mechanics'',
Addison Wesley,
Third edition, 2001
 
\bibitem{Hestenes} David Hestenes,
{\it New foundation for classical mechanics},
Kluwer Academic Publishers,
Second edition, 2003

\bibitem{LandauLifshitz}
L.~D.~Landau and E.~M.~Lifshitz,
``Mecanics'',
Butterworth-Heinenann,
Third edition, 1976 
 
\bibitem{M-S} Mart\'inez R. Sim\'o C.,  {\it On the stability of the Lagrangian homographic solutions in a  curved three body problem on $\mathbb{S}^2$.} Discrete Cont. Dyn. Syst. Ser. A. {\bf 33} (2013), 1157--1175. 

\bibitem{EPC1} P\'erez-Chavela E. and Reyes-Victoria J.G.,
 \emph{An intrinsec approach in the curved $n$-body problem. The positive curvature case}, Trans. Amer. Math. Soc. {\bf 364}-7, (2012),  3805-3827.

\bibitem{tibboel} Tibboel P.,  {\it Polygonal homographic orbits in spaces of constant curvature.} Proc. Amer. Math. Soc.  \textbf{141} (2013), 1465-1471.

\bibitem{Wintner} Wintner A., {\it The
Analytical Foundations Celestial of Mechanics}, Princeton
University Press, Princeton, New York, 1941.
\end{thebibliography}
\end{document}